%&amstex
\input amstex
\NoBlackBoxes
\magnification=\magstep1
 \documentstyle{amsppt}
 %\NoRunningHeads
\rightheadtext{ \\ Explicit Classification for Elliptic Curves }
 \topmatter
 \title
Explicit Classification for Torsion Subgroups\\
of Rational Points of Elliptic Curves
 \endtitle

 \thanks { } \endthanks
 \endtopmatter
 \document

 \hskip 3 cm  Derong QIU \quad  and \quad  Xianke ZHANG
 \par  \smallskip
 \hskip 2cm   (Department of Mathematics, Tsinghua University,
 \par     \smallskip
 \hskip 3.8cm Beijing 100084, P. R. China)
 \par

\vskip 0.6cm

{\bf Abstract.} The classification of  elliptic curves $E$ over the
rationals $\Bbb Q$ is studied according to their torsion subgroups
$E_{tors}(\Bbb Q)$ of rational points. Explicit criteria for
the classification are given when $E_{tors}(\Bbb Q)$
 are  cyclic groups with even orders. The generator points $P$
of $ E_{tors}(\Bbb Q) $  are also explicitly presented
in each case. These results, together with recent results
 of K. Ono, completely solve the problem of the mentioned explicit
 classification when $E$ has a rational point of order 2.
 \vskip 0.06cm
 \



 \qquad arithmetic algebraic geometry \par \vskip 0.8cm

 \heading{I.\quad Introduction and Main Results}\endheading
 \

Let $E$ be an elliptic curve defined over the rationals $\Bbb Q$ .
We consider the following two problems  here:
How to determine
the rational torsion subgroup  from the equation of $E$ ?
How to obtain the generator
 of the rational torsion subgroup explicitly?

From the Mordell-Weil theorem , we know that the Mordell-Weil group
$E(\Bbb Q) $ is a finitely generated abelian group having the form
$$E(\Bbb Q) \cong E_{tors}(\Bbb Q)\times {\Bbb Z}^r $$
where $E(\Bbb Q) $ is the group of the $\Bbb Q$-rational points of $E$,
$E_{tors}(\Bbb Q)$ is the torsion subgroup of $E(\Bbb Q) \ (i.e.,$
all points of $E(\Bbb Q) $  with finite order), $\Bbb Z$ the rational
intergers.

In 1977, B. Mazur completely determined all the possible types of the
rational torsion group
$E_{tors}(\Bbb Q)$ in [1-2]. He shew that the torsion subgroup
must be one of the following fifteen groups:

$${\Bbb Z}/N{\Bbb Z} \qquad (1\le N\le 10 \quad {\text or }\quad N=12);$$
$$ {\Bbb Z}/2{\Bbb Z} \times {\Bbb Z}/2N{\Bbb Z} \qquad (1\le N\le 4).$$
And each of these groups does occur as an $E_{tors}(\Bbb Q)$ for some $E$
(see [4],  p223).

Recently, K.Ono in [3] studied the first  problem above in
an aspect. When torsion subgroups $E_{tors}(\Bbb Q) $ are not cyclic,
he gave criteria to classify them. In fact, K. Ono
considered  elliptic curves $E$ with
$E_{tors}(\Bbb Q)  \supset  {\Bbb Z}/2{\Bbb Z}\times {\Bbb Z}/2{\Bbb Z}$,
 which could be assumed to have the equation
 $E:\quad y^2=x(x+M)(x+N)$, with $M,\ N \in {\Bbb Z}.$
 He obtained necessary and sufficient conditions on $M,\ N$ for the
torsion subgroup to be each of the following types respectively:

\quad (1) $E_{tors}(\Bbb Q)  \supset
{\Bbb Z}/2{\Bbb Z}\times {\Bbb Z}/4{\Bbb Z}$ ;\quad
(2) $E_{tors}(\Bbb Q)  \supset
{\Bbb Z}/2{\Bbb Z}\times {\Bbb Z}/8{\Bbb Z}$  ;

\quad (3) $E_{tors}(\Bbb Q)  \supset
{\Bbb Z}/2{\Bbb Z}\times {\Bbb Z}/6{\Bbb Z}$  .
\par \smallskip

We here will consider the case that the torsion subgroups
$E_{tors}(\Bbb Q)$
 are cyclic with even order,
and will solve the two problems mentioned in the
beginning.

Suppose that
$$E: \quad y^2=f(x)$$
is an elliptic curve, $f(x)\in {\Bbb Q}[x]$. Assume that
$f(x)$ has three complex roots $\alpha,\ \beta,\ \gamma $,
then $P_1=(\alpha , 0),\ P_2=(\beta ,\ 0),\, P_3=(\gamma ,\ 0)$
 are just the three non-trival points of order two of $E$. So
 $E[2]=\{O, P_1, P_2, P_3\}$
 are the group of torsion points of order two of $E$.
By Mazur theorem, $E_{tors}(\Bbb Q) $
is not cyclic if and only if $E_{tors}(\Bbb Q) \supset E[2];\
\ i.e., \ \alpha ,\ \beta ,\ \gamma \in {\Bbb Q}$.
From this we deduce that
$E_{tors}(\Bbb Q) $ is cyclic if and only if
$f(x)$ has at most one rational root.
Thus we know that  $E_{tors}(\Bbb Q) $ is cyclic with even order
if and only if $E$ has just one non-trival rational
 point of order 2; that is,
 $f(x)$ has just one rational root. Let us assume so. Then, up to
a translation, we may assume this rational root
of $f(x)$ to be $0$; so we have
$f(x)=x(x-\alpha )(x-\beta )  $. Since $\alpha +\beta $ and
$\alpha \beta $
are rationals, we have  $\alpha = a+b \sqrt D, \ \beta=a-b\sqrt D,$
where $ a,b\in {\Bbb Q},\ b\not= 0,\ D $ a squarefree integer.
 We may further assume $a, b$ to be rational integers and
  gcd$\{a,b\}=(a,b)$ is squarefree. In fact, when we replace $x,\ y $
by $x/d^2,\ y/d^3$ , then the equation of
$E:\ y^2=x(x-\alpha )(x-\beta)$ becomes
$E_d:\ y^2=x(x-d^2\alpha )(x-d^2\beta). \ E$ and $E_d$ are
${\Bbb Q}$-isomorphic, so
$E(\Bbb Q)\cong E_d(\Bbb Q),\,
E_{tors}(\Bbb Q) \cong {E_d}_{tors} (\Bbb Q)$.
 Our main result here is the following theorem.

\vskip 0.25cm
\quad \proclaim { \quad THEOREM 1}   Suppose that
$E:\ y^2=x(x+M)(x+N)$
is an elliptic curve, where $M=m+n\sqrt D,\ N=m-n\sqrt D,\; D $ and
$(m,\ n)$ are squarefree integers, $\; D\not=1,\ n\not=0,$
 and $m$ are all rational integers. Then the $\Bbb Q$-rational torsion
subgroup $E_{tors}(\Bbb Q) $ of $E$ is classified as follows:
\smallskip

(I) $E_{tors}(\Bbb Q) \supset \Bbb Z /4\Bbb Z $ if and only if
$m=a^2+b^2D,\ n=2ab,\ $ where $a,b\in \Bbb Z $
are relatively prime and non-zero.
                 \smallskip

(II) $E_{tors}(\Bbb Q) = \Bbb Z /8\Bbb Z $ if and only if
$m=u^4+v^2 w^2 D,\ n=2u^2vw,\ 2u^2-v^2=w^2D $,
 where $u,v,w\in \Bbb Z $ are non-zero.
                 \smallskip

(III) $E_{tors}(\Bbb Q) \supset \Bbb Z /6\Bbb Z $ if and only if
$m=a^2+2ac+b^2D,\ n=2b(a+c),\ a^2-b^2D=c^2$,
 where $a,b,c\in\ \Bbb Z $ are relatively prime and non-zero.
                 \smallskip

(IV) $E_{tors}(\Bbb Q) = \Bbb Z /12\Bbb Z $ if and only if
$m=v^2-u^2+ w^2 D,\ n=2vw,\   $ and \newline
$ 3(v^2-w^2D)^4-4u^2(v^2-w^2D)^2(v^2+w^2D)-16u^4v^2w^2D=0$,
 where $u,v,w\in \Bbb Z $ \newline
are non-zero.
                 \smallskip

(V) $E_{tors}(\Bbb Q) = \Bbb Z /10\Bbb Z $ if and only if
$m=2s(s+u)-v^2,\ n=2st,\ ( s+u)^2-v^2=t^2D, $  and
$(u-v)^2(u+v)=4uvs,\ $
 where $u,v,s,t\in \Bbb Z $ are non-zero.
                 \smallskip

(VI) Otherwise, $E_{tors}(\Bbb Q) = \Bbb Z /2\Bbb Z $ .  \endproclaim

\vskip 0.25cm
\quad \proclaim {\quad THEOREM 2}   Let  $P_n$
denote a generator (point) of the torsion group $E_{tors}(\Bbb Q) $
of the rational points of order $n$ in $E$.
Then $x(P_n)$ and $x(2P_n)$  , the $x$-coordinates of $P_n$ and $2P_n$
respectively,  could be as
in the following, where the  cases and notations are corresponding to
Theorem 1:
\smallskip

(I) $x(P_4)=a^2-b^2D;\quad x(2P_4)=0.$\smallskip

(II) $x(P_8)=(u+v)(v-u)^3;\quad x(2P_8)=(u^2-v^2)^2.$\smallskip

(III) $x(P_6)=5c^2+4ac;\quad x(2P_6)=c^2.$\smallskip

(IV)  $x(P_{12})=(u+v)^2-w^2D;\quad x(2P_{12})=u^2.$\smallskip

(V) $x(P_{10})=2v^2+4vs-u^2;\quad  x(2P_{10})=u^2.$\smallskip

(VI) $x(P_2)=0.$     \endproclaim
\vskip 0.15 cm

\heading  {II.\hskip 0.1cm Proofs of the Theorems} \endheading

\quad \proclaim {\quad Lemma 1}   Suppose that
$E:\ y^2=(x-\alpha )(x-\beta)(x-\gamma)$ is an elliptic curve over
any number field $K$,  $\alpha ,\beta ,\gamma \in K$.
Let the point $(x_0,\ y_0)\in E(K)$. Then $(x_0,\ y_0)=2(x_1,\ y_1)$
for some point $(x_1,\ y_1)\in E(K)$ if and only if
$x_0-\alpha ,\ x_0-\beta ,\ $ and $x_0-\gamma \ $
all are squares in $K$
(see [5], p85).
\par \smallskip \endproclaim

\quad \proclaim  { \quad Proof of Theorem 1}  By Lutz-Nagell
theorem,  any rational torsion point
 $P\in E_{tors}(\Bbb Q) $ is an integer point, $i.e., $ the coordinates
$x(P),\  y(P) \in \Bbb Z $
(see [5]). The following
duplication formula could be obtained from  formulae of [4]:
$$x(2P)=((x(P)^2-MN)/2y(P))^2 .  \eqno(*)$$

(I) If $E_{tors}(\Bbb Q) \supset \Bbb Z /4\Bbb Z $ , then $E$ has a rational
point $P$ of order 4,
and $2P=P_0=(0,0)$ is the unigue rational point of order 2.
 By Lemma 1,   $ M$ and $ N$
 are squares in the field $K=\Bbb Q (\sqrt D)$.
In fact we could easily see that
$M,\ N$ are squares in the ring $\Bbb Z [\sqrt D]$.
In other words, we have
$M=(a+b\sqrt D)^2,\ N=(a-b\sqrt D)^2,\; a,b\in \Bbb Z$.
So we have $m=a^2+b^2D,\; n=2ab$. Since $(m,n)$ is squarefree and
$n\not= 0$, so $(a,b)=1,\; ab \not= 0$.

Conversely, if the conditions on $m,n$ in (I) hold,
 then $M,N$ are squares in $K$.
By Lemma 1, there is a $K$-rational point $P$ of $E$ such that
 $2P=P_0=(0,0)$. By the duplication formula (*) we have
$x(P)^2-MN=0,\;\ x(P)^2=MN=(a^2-b^2D)^2,\
 x(P)=\pm (a^2-b^2D).$ Substitute this
$x(P)$ into the equation of $E$, we obtain two integer points $P$:
$x(P)= (a^2-b^2D),\ y(P)=\pm 2a(a^2-b^2D)\ ( $
There is no rational point $P$
with $x(P)=-(a^2-b^2D))$ .  These  $P$ are points of order 4 in
  $E_{tors}(\Bbb Q) $.
\par \smallskip

(II) Suppose that $E_{tors}(\Bbb Q) = \Bbb Z /8\Bbb Z $ and
 $P$ is a rational point of order 8
of $E$. So $2P$ is of order 4, and by (I) we know that
$m=a^2+b^2D,\ n=2ab,\ $  and $x(2P)=a^2-b^2D$.
Then by Lemma 1 we know that
$x(2P)=a^2-b^2D,\  x(2P)+M=2a^2+2ab\sqrt D ,\  x(2P)+N=
2a^2-2ab\sqrt D $ all are squares in field $K$. From the
duplication formula (*) and the fact that
 $x(P),\ y(P),\  x(2P),\ MN,\ a, $ and $b$
are rational integers, we obtain that
(i) $a^2-b^2D=c^2\  $;  (ii) $2a^2+2ab\sqrt D=(s+t\sqrt D)^2$; where
$c,\ s,\ t \in \Bbb Z$. These give: (iii) $s^2=a(a+c)$;
(iv) $t^2D=a(a-c)$. Since $(a,b)=1$ and $D$ is squarefree,
 from (i) we have $(a,c)=1$. So $(a, a+c)=1,$ and via (iii) we have
$a=u^2,\ a+c=v^2$, where $u,v\in \Bbb Z $ and $(u,v)=1$. Note that
$D$ is squarefree, so by (iv) we have
$2u^2-v^2=w^2D $ with $w\in \Bbb Z $.
Then via (i) we have $b=vw$,
so $m=u^4+v^2w^2D,\ n=2u^2vw,\ 2u^2-v^2=w^2D,$ where
$u,v,w\in \Bbb Z$ are non-zero.

   Conversely, suppose that $E$ satisfies the given condition.
 So $E$ also satisfies the condition in Case (I). Thus $E$ contains a
 $\Bbb Q $-rational point $P_4$ of order 4 with
$x(P_4)=(u^2-v^2)^2$. It is easy to verify
that the coordinates of $P_4$ satisfy Lemma 1 for
$K=\Bbb Q (\sqrt D)$.
So there is a $K$-rational point $P$ with $2P=P_4$, and then $P$ has
order 8 and $x(2P)=x(P_4)=(u^2-v^2)^2$. From the duplication formula (*)
we have $4y^2(u^2-v^2)^2=(x^2-MN)^2$  (here $x=x(P)$).
Substitute the relations for $m,n$  into the this equation and the
equation of $E$ we obtain
$$ x^4-4(u^2-v^2)^2x^3-2(u^2-v^2)^2(5u^4+6u^2v^2-3v^4)x^2-
4(u^2-v^2)^6x+(u^2-v^2)^8=0,$$
which turns to be
$$(x-(u^2-v^2)^2)^4=16u^4(u^2-v^2)^2x^2, $$
having a rational-integer solution
$$x=(u+v)(v-u)^3. $$
It is easy to see that this is the only integer solution.
Substituting this $x$ into the equation of $E$, we obtain $y$ which is
obviously a rational-integer. So we find a
$\Bbb Q$-rational point $P_8$ of order 8, and obtain that
 $x(P_8)=(u+v)(v-u)^3$. This proves case (II).
\par \smallskip

(III) Suppose that $E_{tors}(\Bbb Q) \supset \Bbb Z /6\Bbb Z $ .
 Then there is a rational
point of order 3 in $E$, denote it by $P$.
Obviously $x(2P)=x(P)\not =0$.
From the duplication formula (*) we have
$x(2P)=u^2=x(P),\ u\in \Bbb Z $.
By the property of torsion points with order 3
(see [6], p40), we know that x=x(P) satisfies the equation
$3x^4+4(M+N)x^3+6MNx^2-M^2N^2=0$.
As a homogeneous polynomial equation of degree 4 in the
variables $M,\ N, \ and x$, it could be parametrized
(due to Nigel Boston, see [3]) as
$ M/x=(1+t)^2-1,\ N/x=(1+1/t)^2-1$ for some $t$.
for some $t$. Obviously $t\in \Bbb Q (\sqrt D)-\Bbb Q,$ Let
$t=(a+b\sqrt D)/c,\ $ with $(a,b,c)=1,\ bc\not=0, \ a,b,c\in \Bbb Z $.
Substituting $t$ and the above $M/x,\  N/x$ into the equation of $E$,
we have $(2+a/c+ac/(a^2-b^2D))(b/c-bc/(a^2-b^2D)=0$.
If $2+a/c+ac/(a^2-b^2D)=0$, then $2+t+1/t=(b/c-bc/(a^2-b^2D)\sqrt D$.
Substituting into the equation of $E$, we have
$(y/x)^2/x=(b/c-bc/(a^2-b^2D)^2D$, which is
impossibe since $x=u^2 $ and $D$ is squarefree.
Therefore we must have $b/c-bc/(a^2-b^2D)=0, \  a^2-b^2D=c^2$. So
$1/t = (a-b\sqrt D)/c,$ and we obtain
$$m=x(a^2+2ac+b^2D)/c^2,\quad n=2xb(a+c)/c^2.$$
It is easy to verify that $(a^2+2ac+b^2D,\ 2b(a+c),\ c^2)$
is squarefree, which implies $x=c^2$. So
$m=a^2+2ac+b^2D,\  n=2b(a+c),\ a^2-b^2D=c^2,$ as desired.

  Now suppose that $E$ satisfies the given condition of (III).
 From the condition we could easily obtain a
 $\Bbb Q -$rational point $P_3$ of order 3 with
$x(P_3)=c^2$ and $|y(P_3)|=2|a+c|c^2$ (Actually, every rational point
of order 3 of such $E$ satisfies $x(P_3)=c^2$ ). Thus
$P_3+P_0=P_6$ is a non-trival $\Bbb Q -$rational point of order 6,
where $P_0=(0,0)$ is a point of order 2. Then via the coordinate formula
for the group law of an elliptic curve, we obtain
$x(P_6)=5c^2+4ac$. This proves the case (III).
 \par \smallskip

(IV)Suppose that $E_{tors}(\Bbb Q) = \Bbb Z /12\Bbb Z $ .
Then $E$ has a rational point $P$
of order 12. So $2P$ has order 6.  Similarly as we just proved case (II)
by using some results of case (I), we could obtain the following by using
results of case (III):     $\ m=v^2-u^2+w^2D,\ n=2vw$. We thus also have
(1) $\ b(a+c)=vw; \ $ (2) $5c^2+4ac=u^2 ;\ $ (3) $ a^2-b^2D=c^2;\ $ and
(4) $4c(a+c)+(a+c)^2+b^2D=v^2+w^2D$; where
 $a,b,c $ are as in case (III), and
$u,v,w\in \Bbb Z $ are non-zero. From these formulae (1--4) and via
calculation, we could obtain the desired expression
$$3(v^2-w^2D)^4-4u^2(v^2-w^2D)^2(v^2+w^2D)-16u^4v^2w^2D=0.$$

Conversely, suppose that $E$ satisfies the given condition
in (IV). Then $E$
satisfies the condition of case (III), so $E$ has a rational point
$P_6$ of order 6 with $ x(P_6)=5c^2+4ac=u^2$.
By a method similar to that in case (II),
we  obtain that $E$ has a $\Bbb Q (\sqrt D)-$rational point
$P$ of order 12  and $x(P)$ satisfies the following equation:
$$x^4-4u^2x^3-hx^2-4u^2ex+e^2=0$$
where
$$h=2((v^2-w^2D)^2+2u^2(v^2+w^2D)-3u^4),$$
$$e=(v^2-w^2D)^2+u^4-2u^2(v^2+w^2D).$$
Via a careful calculation using the condition satisfied by $E$, we
deduce the above equation as
$$((x-u^2)^2-(2u^2(v^2+w^2D)-(v^2-w^2D)^2))^2=4(v^2-w^2D)^2x^2.$$
Thus we obtain the integer solution $x=(u+v)^2-w^2D$, and then get
$y=y(P)\in \Bbb Z$. So $E$ has a non-trival $\Bbb Q-$rational
 point $P_{12}$ of order 12, and
$x(P_{12})=(u+v)^2-w^2D$. The case (IV) is proved.
\par \smallskip

(V) Suppose that $E_{tors}(\Bbb Q) = \Bbb Z /10\Bbb Z $
and $P$ is a non-trival rational
point of order 5. Then it is easy to see that $x(4P)=x(P)$ and
$x(2P)\not=x(P)$. Then using the duplication formula (*) as in the
formal cases we could obtain that $m=2s(s+u)-v^2,\ n=2st $, and
$(s+u)^2-v^2=t^2D,\ (u+v)(u-v)^2=4uvs,$ where $u,v,s,t\in \Bbb Z $
are non-zero.

Conversely, suppose that $E$ satisfies the condition of (V). Then from
the condition it is easy to find a rational point $P_5$ of order 5 and
obtain $x(P_5)=u^2,\ |y(P_5)|=|u(u^2-v^2+2us)|$. So $P_0+P_5=P_{10}$ is
a rational point of order 10. Then via the coordinate formulae of
the additive law of elliptic curves, we obtain
$x(P_{10})=2v^2+4vs-u^2$.
This completes the proof of Theorem 1.\qed \endproclaim

In the above proof for Theorem 1, we have also proved
 the results of theorem 2.

 \par
 \newpage
 \
 \heading References  \phantom {xxxxxxxxxxxxxxxxxxxxxxxxxxxxxxxxxxxxxxxxxxxxxx
xxxxxxxxxxxxxxxxxxx} \endheading

 [1]  B\. Mazur, Modular curves and the Eisenstein ideal,
IHES Publ. Math. 47(1977), 33-186 .\newline
  \

[2] B\. Mazur, Rational points on modular curves, Modular
Functions of One Variable V, Lecture Notes in Math. 601(1977),
 107-148, Springer-Verlag, New York.\newline
\

 [3] K\. Ono, Euler's concordant forms, Acta Arithmetica,
LXX VIII 2(1996), 101-123. \newline
 \

 [4] J. Silverman, The Arithmetic of Elliptic Curves, Springer-Verlag,
New York, 1986. \newline
 \

 [5] A. knapp, Elliptic Curves, Princeton Univ. Press,
Princeton, 1992.\newline
 \

 [6] J. Silverman and J.Tate, Rational Points on Elliptic
Curves, Springer-Verlag, New York, 1992. \newline

 \

 \medskip

 Department of Mathematics, Tsinghua University, Beijing 100084,
P. R. China

E-mail: xianke\@tsinghua.edu.cn

\hskip 1.3cm  xzhang\@math.tsinghua.edu.cn

Fax:  0086-10-62562768

(This paper was published in:

Acta Mathematica. Sinica (English Series), 18(2002.7), No.3,
539-548)

 \enddocument

\end